\documentclass[a4paper]{article}

\usepackage[utf8]{inputenc}
\usepackage{lmodern}
\usepackage{amssymb,amsfonts,amsmath,amsthm,bbm,mathrsfs,aliascnt,mathtools,pifont}
\usepackage[english]{babel}
\usepackage[colorlinks]{hyperref}
\usepackage{tikz}
\usetikzlibrary{decorations}
\usetikzlibrary{decorations.pathreplacing}
\tikzset{individu/.style={draw,thick}}

\usepackage{pgfplots}

\theoremstyle{plain}
\newtheorem{theorem}{Theorem}[section]

\newtheorem{lemma}[theorem]{Lemma}
\newtheorem{proposition}[theorem]{Proposition}

\theoremstyle{definition}

\theoremstyle{remark}
\newtheorem{remark}[theorem]{Remark}
\newtheorem{example}[theorem]{Example}

\numberwithin{equation}{section}

\newcommand{\N}{\mathbb{N}}

\newcommand{\R}{\mathbb{R}}
\newcommand{\C}{\mathbb{C}}

\newcommand{\X}{\mathbf{X}}

\newcommand{\indset}[1]{\mathbf{1}_{#1}}

\newcommand{\dd}{\mathrm{d}}

\DeclareMathOperator{\E}{\mathbb{E}}

\renewcommand{\P}{\mathbb{P}}

\renewcommand{\rho}{\varrho}
\renewcommand{\epsilon}{\varepsilon}

\newcommand\Zz{\mathbf{Z}}

\newcommand\Pp{\mathcal{P}}
\newcommand\Mm{\mathcal{M}}

\newcommand\m{{\boldsymbol{m}}}
\newcommand\s{{\mathbf s}}
\newcommand\Ss{{\mathbf S}}

\makeatletter \newcommand\mathof[1]{{\operator@font#1}} \makeatother

\begin{document}

\title{On a population model
with memory }
\author{Jean Bertoin\footnote{Institute of Mathematics, University of Zurich, Switzerland, \texttt{jean.bertoin@math.uzh.ch}}  }
 \date{\small}
\maketitle 
\thispagestyle{empty}

\begin{abstract} 
{Consider first a memoryless population model described  by the usual branching process  with a given mean reproduction matrix on a finite space of types.
Motivated by the consequences of atavism in Evolutionary Biology, we are interested in a modification of the dynamics where  individuals keep full memory of their forebears and procreation  involves the reactivation of a gene picked at random on the ancestral lineage.
By comparing the spectral radii of the two mean reproduction matrices (with and without memory), we observe that, on average, the model with memory always grows at least as fast as the model without memory. The proof relies on analyzing a biased Markov chain on the space of memories, and the existence of a unique ergodic law is demonstrated through asymptotic coupling.
}
  \vskip 1mm
{\normalfont \bfseries Keywords}: Multitype branching process, Markov chains with memory, spectral radius, atavism.\newline
\vskip 1mm
{\normalfont \bfseries Mathematics Subject Classification:}  60J80; 60F99

\end{abstract}

\section{Introduction}\label{s:intro}

The genetic material that is transmitted  in living organisms from parents to children through DNA keeps memory of past events. It is known that traits that have disappeared phenotypically do not necessarily disappear from DNA and some genes can remain dormant for a great many generations. 
Atavism, also known as evolutionary throwbacks, or reversions, is indeed a well-known phenomenon in Evolutionary Biology, which refers to the reappearance of a trait that was last expressed in remote ancestors\footnote{A precise definition of atavism given in \cite{Wagner}
 reads ``a morphophysiological trait controlled by an ancient gene regulatory network that survives periods of disuse and can be reactivated and reused in later lineages, even if it was not used in the immediate ancestors" .}.  This can typically occur when inactive genes preserved in the genome become activated, for instance by a mutation or by a fault in the genetic control suppressing the expression of that gene. Broadly speaking, a common view in this area is that atavism enables to look back on evolutionary solutions that were once beneficial and could become useful again, and therefore strengthens the adaptability and survival potential of the population.

Our purpose here is to point out on a simple population model, that even in the absence of natural selection and  environmental changes which might favor old traits again,
 the sole possibility of blindly re-activating genetic material from ancestors is always beneficial in comparison with the model without memory.  
 In other words, the advantage of atavism for a population does not only stem from improving its adaptability, there is always an intrinsic benefit induced by the sole use of the memory from preceding generations for reproduction. Heuristically, this can be surmised as genetic material from prolific remote ancestors remains present in a larger part of the current  population and is thus more likely to be re-activated, producing in turn a larger descent in the future. 

The benchmark population model without memory is the classical multitype branching process,  where the space of types is a finite set  $S$.
Its evolution is governed by the basic reproduction kernel  
$$\pi(s, \cdot), \qquad s\in S,$$ 
that specifies the distribution of the types in the progeny given the type of the parent. 
That is, each $\pi(s, \cdot)$ is a  probability measure on  $\N^S$, and for any vector with integer coordinates  $n=(n_t)_{t\in S}$ in $\N^S$,
the probability that a parent with type $s$ begets $n_t$ children of type $t$ for every $t\in S$ equals $\pi(s, n)$. 
It is well-known that spectral properties of the mean reproduction matrix $m\in \R_+^{S\times S}$ with entries
\begin{equation} \label{E:meanm}
m(s,t) = \sum_{n\in \N^S} n_t \pi(s,t), \qquad s,t\in S
\end{equation}
have a key role for the large time behavior of the branching process; 
we assume throughout this work and without further mentions that $m$ is irreducible and aperiodic.
Then the growth rate of the population can be identified 
as the spectral radius $r$ of $m$;  $r$ is also know as the principal (or Perron-Frobenius) eigenvalue. Furthermore the probability vector on $S$ which describes   the asymptotic proportion of types is an eigenvector for the principal eigenvalue of the transposed mean reproduction matrix. 
 See, for instance \cite[Chapter V]{AN}, or \cite[Chapter II]{Harris}.

We introduce now a simple variation for the evolution of populations where individuals keep  the memory of the types of their forebears and then activate a randomly picked type for their reproduction. This bears some similarities with so-called reinforced Galton--Watson processes, which have been recently introduced in \cite{BM1}, see also \cite{BM2}.
Specifically,  imagine that every individual carries an infinite sequence 
$\s=(s_0, s_1, \ldots)$, which we think of as its memory, where each  $s_j$ is the type of its forebear $j$ generations backwards, starting with the type $s_0$ of the individual itself.  The offspring of an individual with memory $\s$ can then be represented by a vector $n=(n_t)_{t\in S}$ in $\N^S$, meaning that for every $t\in S$, 
  $n_t$ children are born with memory simply obtained by concatenating the  type $t$ of the child with the memory $\s$ of the parent, i.e. $t\s=(t, s_0, s_1, \ldots )$.
In addition to the basic reproduction kernel  $\pi$ which was already introduced above, 
the second item needed to specify the evolution of the population with memory is a probability measure on $\N$,
$$\tau=(\tau(j): {j\geq 0}),$$ 
that describes the statistics of the random (backward) generation at which the memory is activated for procreation.
We always assume that $\tau(0)>0$ and that the first moment of $\tau$ is finite,
$$\sum_{j=0}^{\infty} j \tau(j)<\infty.$$
We mention that when the support of $\tau$, $\mathrm{Supp}(\tau)=\{j\geq 0: \tau(j)>0\}$, is finite, it suffices actually to consider memory sequences with length $\max \mathrm{Supp}(\tau)$, and several technical aspects of our work can then be made simpler. However, it is a known fact in Evolutionary Biology that ancestral traits may reappear after being dormant for millions of years\footnote{As it is shown for instance by the existence of extremely rare cases of ``true human tail'',
or the formation of reptilian-like patterns such as teeth in chicken after experimental manipulation.},  so 
it seems more natural to include in our study cases where $\mathrm{Supp}(\tau)$ is unbounded. Further, solving the technical difficulties that unboundedness of the memory induces may also be interesting mathematically in its own  right.

Thus,  when an individual with memory $\s$ procreates,  a random type $s_T$ from its memory  is activated, where $T$ is a random  variable distributed according to $\tau$
(so the model without memory is recovered by choosing $\tau=\delta_0$, i.e. $T\equiv 0$ a.s.). Then, conditionally on $s_T=s$, the law of the random vector  in $\N^S$ that records the types of its progeny is $\pi(s,\cdot)$. 
The counterpart  in the case with memory of the mean reproduction matrix $m$ is given by the operator $\boldsymbol{m}$ acting on the space of bounded measurable functions $f: S^\N \to \R$,
   \begin{equation} \label{E:T1}
   \boldsymbol{m}f(\s)= \sum_{j=0}^{\infty} \tau(j) \sum_{t\in S} m(s_j,t) f(t\s), \qquad \s \in S^\N.
   \end{equation}
We also assume that independent copies of $T$ are used for different individuals, 
and  that the branching property is fulfilled. In other words, conditionally on the genetic memories of their parents at some generation, the memories of children at the next generation are independent. So the population model with memory can be viewed as a branching process on the space  $S^\N$ of sequences of types.

 We stress that, roughly speaking, the activation of the memory is fair, in the sense that the law of $T$ does not depend on the memory, and hence does not favor \textit{a priori} certain types more than others (which could be the case if the law of $T$ depended on the memory of the parent, allowing for instance the activation of an advantageous type in the memory).
  The central result of this work
states that the spectral radius $\boldsymbol{r}$ of the operator $\boldsymbol{m}$ is always bounded from below by the  spectral radius $r$ of the mean reproduction matrix $m$.
Since average growth rates of general branching processes are determined by the spectral radii of their mean reproduction operators, this shows that the model with memory always growths on average at least  as fast as the model without memory, hence stressing that remembering is always beneficial. 
\begin{theorem} \label{T:1}
There is the inequality
$\boldsymbol{r}\geq r$.
\end{theorem}
We will also point at a more elementary upper-bound for $\boldsymbol{r}$, see the forthcoming \eqref{E:Harn},
and at situations where the inequality of Theorem \ref{T:1} is actually an equality, see the forthcoming Example \ref{E:bal}.
On the other hand, simple examples as the one below for which explicit calculations can be performed show that  the inequality is strict in general (and possibly rather sharp). 
\begin{example} Consider the case $S=\{a,b\}$, $m=\left(\begin{matrix}1 & 1\\ 1 & 2\\ \end{matrix}\right)$, and $\tau(0)=u$, $\tau(1)=1-u$ for some parameter  $u\in(0,1)$.
If we enumerate the four memory sequences of length 2 in the lexicographic order,  $(a,a), (a,b), (b,a)$ and $(b,b)$, then we get
$$\boldsymbol{m}=\boldsymbol{m}(u)=\left(\begin{matrix}1 & 0 & 1 & 0\\ 1 &0 & 2-u &0\\
0 & 1 & 0 & 1+u\\
0 & 1 & 0 & 2\\ \end{matrix}\right). $$
The spectral radius of $m$ is $r=(3+\sqrt 5)/2 \approx 2.618$; 
the plot of the spectral radius $\boldsymbol{r}(u)$ of $\boldsymbol{m}(u)$ as a function of $u$ is displayed  below.
\vskip 6mm
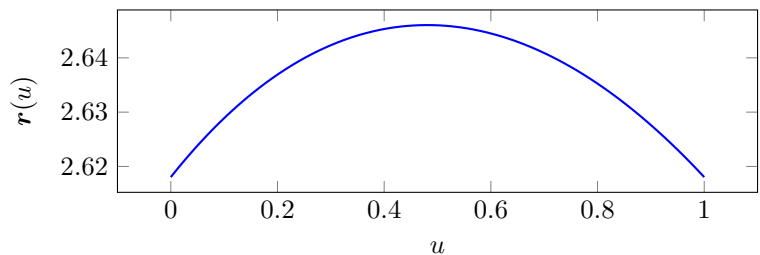
\begin{tikzpicture}
    \begin{axis}[
        xlabel={$u$},
        ylabel={$\boldsymbol{r}(u)$},
        grid=minor,
        width=10cm,
        height=4cm,
    ]
        \addplot[smooth, thick, blue] coordinates {
(0.00000, 2.61803) (0.00200, 2.61828) (0.00401, 2.61853) (0.00601, 2.61877) (0.00802, 2.61901) (0.01002, 2.61926) (0.01202, 2.61950) (0.01403, 2.61974) (0.01603, 2.61998) (0.01804, 2.62021) (0.02004, 2.62045) (0.02204, 2.62068) (0.02405, 2.62092) (0.02605, 2.62115) (0.02806, 2.62138) (0.03006, 2.62161) (0.03206, 2.62184) (0.03407, 2.62207) (0.03607, 2.62229) (0.03808, 2.62252) (0.04008, 2.62274) (0.04208, 2.62296) (0.04409, 2.62319) (0.04609, 2.62341) (0.04810, 2.62362) (0.05010, 2.62384) (0.05210, 2.62406) (0.05411, 2.62427) (0.05611, 2.62449) (0.05812, 2.62470) (0.06012, 2.62491) (0.06212, 2.62512) (0.06413, 2.62533) (0.06613, 2.62554) (0.06814, 2.62575) (0.07014, 2.62596) (0.07214, 2.62616) (0.07415, 2.62636) (0.07615, 2.62657) (0.07816, 2.62677) (0.08016, 2.62697) (0.08216, 2.62717) (0.08417, 2.62737) (0.08617, 2.62756) (0.08818, 2.62776) (0.09018, 2.62795) (0.09218, 2.62815) (0.09419, 2.62834) (0.09619, 2.62853) (0.09820, 2.62872) (0.10020, 2.62891) (0.10220, 2.62909) (0.10421, 2.62928) (0.10621, 2.62947) (0.10822, 2.62965) (0.11022, 2.62983) (0.11222, 2.63002) (0.11423, 2.63020) (0.11623, 2.63038) (0.11824, 2.63055) (0.12024, 2.63073) (0.12224, 2.63091) (0.12425, 2.63108) (0.12625, 2.63126) (0.12826, 2.63143) (0.13026, 2.63160) (0.13226, 2.63177) (0.13427, 2.63194) (0.13627, 2.63211) (0.13828, 2.63228) (0.14028, 2.63245) (0.14228, 2.63261) (0.14429, 2.63278) (0.14629, 2.63294) (0.14830, 2.63310) (0.15030, 2.63326) (0.15230, 2.63342) (0.15431, 2.63358) (0.15631, 2.63374) (0.15832, 2.63389) (0.16032, 2.63405) (0.16232, 2.63420) (0.16433, 2.63436) (0.16633, 2.63451) (0.16834, 2.63466) (0.17034, 2.63481) (0.17234, 2.63496) (0.17435, 2.63511) (0.17635, 2.63526) (0.17836, 2.63540) (0.18036, 2.63555) (0.18236, 2.63569) (0.18437, 2.63583) (0.18637, 2.63597) (0.18838, 2.63612) (0.19038, 2.63625) (0.19238, 2.63639) (0.19439, 2.63653) (0.19639, 2.63667) (0.19840, 2.63680) (0.20040, 2.63694) (0.20240, 2.63707) (0.20441, 2.63720) (0.20641, 2.63733) (0.20842, 2.63746) (0.21042, 2.63759) (0.21242, 2.63772) (0.21443, 2.63785) (0.21643, 2.63797) (0.21844, 2.63810) (0.22044, 2.63822) (0.22244, 2.63834) (0.22445, 2.63847) (0.22645, 2.63859) (0.22846, 2.63871) (0.23046, 2.63883) (0.23246, 2.63894) (0.23447, 2.63906) (0.23647, 2.63918) (0.23848, 2.63929) (0.24048, 2.63940) (0.24248, 2.63952) (0.24449, 2.63963) (0.24649, 2.63974) (0.24850, 2.63985) (0.25050, 2.63996) (0.25251, 2.64007) (0.25451, 2.64017) (0.25651, 2.64028) (0.25852, 2.64038) (0.26052, 2.64049) (0.26253, 2.64059) (0.26453, 2.64069) (0.26653, 2.64079) (0.26854, 2.64089) (0.27054, 2.64099) (0.27255, 2.64109) (0.27455, 2.64118) (0.27655, 2.64128) (0.27856, 2.64137) (0.28056, 2.64147) (0.28257, 2.64156) (0.28457, 2.64165) (0.28657, 2.64174) (0.28858, 2.64183) (0.29058, 2.64192) (0.29259, 2.64201) (0.29459, 2.64209) (0.29659, 2.64218) (0.29860, 2.64226) (0.30060, 2.64235) (0.30261, 2.64243) (0.30461, 2.64251) (0.30661, 2.64259) (0.30862, 2.64267) (0.31062, 2.64275) (0.31263, 2.64283) (0.31463, 2.64291) (0.31663, 2.64298) (0.31864, 2.64306) (0.32064, 2.64313) (0.32265, 2.64321) (0.32465, 2.64328) (0.32665, 2.64335) (0.32866, 2.64342) (0.33066, 2.64349) (0.33267, 2.64356) (0.33467, 2.64362) (0.33667, 2.64369) (0.33868, 2.64376) (0.34068, 2.64382) (0.34269, 2.64388) (0.34469, 2.64395) (0.34669, 2.64401) (0.34870, 2.64407) (0.35070, 2.64413) (0.35271, 2.64419) (0.35471, 2.64425) (0.35671, 2.64430) (0.35872, 2.64436) (0.36072, 2.64441) (0.36273, 2.64447) (0.36473, 2.64452) (0.36673, 2.64457) (0.36874, 2.64462) (0.37074, 2.64467) (0.37275, 2.64472) (0.37475, 2.64477) (0.37675, 2.64482) (0.37876, 2.64487) (0.38076, 2.64491) (0.38277, 2.64496) (0.38477, 2.64500) (0.38677, 2.64504) (0.38878, 2.64509) (0.39078, 2.64513) (0.39279, 2.64517) (0.39479, 2.64521) (0.39679, 2.64524) (0.39880, 2.64528) (0.40080, 2.64532) (0.40281, 2.64535) (0.40481, 2.64539) (0.40681, 2.64542) (0.40882, 2.64546) (0.41082, 2.64549) (0.41283, 2.64552) (0.41483, 2.64555) (0.41683, 2.64558) (0.41884, 2.64561) (0.42084, 2.64563) (0.42285, 2.64566) (0.42485, 2.64568) (0.42685, 2.64571) (0.42886, 2.64573) (0.43086, 2.64576) (0.43287, 2.64578) (0.43487, 2.64580) (0.43687, 2.64582) (0.43888, 2.64584) (0.44088, 2.64586) (0.44289, 2.64587) (0.44489, 2.64589) (0.44689, 2.64590) (0.44890, 2.64592) (0.45090, 2.64593) (0.45291, 2.64595) (0.45491, 2.64596) (0.45691, 2.64597) (0.45892, 2.64598) (0.46092, 2.64599) (0.46293, 2.64600) (0.46493, 2.64600) (0.46693, 2.64601) (0.46894, 2.64602) (0.47094, 2.64602) (0.47295, 2.64602) (0.47495, 2.64603) (0.47695, 2.64603) (0.47896, 2.64603) (0.48096, 2.64603) (0.48297, 2.64603) (0.48497, 2.64603) (0.48697, 2.64603) (0.48898, 2.64602) (0.49098, 2.64602) (0.49299, 2.64602) (0.49499, 2.64601) (0.49699, 2.64600) (0.49900, 2.64599) (0.50100, 2.64599) (0.50301, 2.64598) (0.50501, 2.64597) (0.50701, 2.64596) (0.50902, 2.64594) (0.51102, 2.64593) (0.51303, 2.64592) (0.51503, 2.64590) (0.51703, 2.64589) (0.51904, 2.64587) (0.52104, 2.64585) (0.52305, 2.64584) (0.52505, 2.64582) (0.52705, 2.64580) (0.52906, 2.64578) (0.53106, 2.64575) (0.53307, 2.64573) (0.53507, 2.64571) (0.53707, 2.64568) (0.53908, 2.64566) (0.54108, 2.64563) (0.54309, 2.64561) (0.54509, 2.64558) (0.54709, 2.64555) (0.54910, 2.64552) (0.55110, 2.64549) (0.55311, 2.64546) (0.55511, 2.64543) (0.55711, 2.64540) (0.55912, 2.64536) (0.56112, 2.64533) (0.56313, 2.64529) (0.56513, 2.64526) (0.56713, 2.64522) (0.56914, 2.64518) (0.57114, 2.64514) (0.57315, 2.64510) (0.57515, 2.64506) (0.57715, 2.64502) (0.57916, 2.64498) (0.58116, 2.64494) (0.58317, 2.64489) (0.58517, 2.64485) (0.58717, 2.64480) (0.58918, 2.64476) (0.59118, 2.64471) (0.59319, 2.64466) (0.59519, 2.64461) (0.59719, 2.64456) (0.59920, 2.64451) (0.60120, 2.64446) (0.60321, 2.64441) (0.60521, 2.64435) (0.60721, 2.64430) (0.60922, 2.64425) (0.61122, 2.64419) (0.61323, 2.64413) (0.61523, 2.64408) (0.61723, 2.64402) (0.61924, 2.64396) (0.62124, 2.64390) (0.62325, 2.64384) (0.62525, 2.64378) (0.62725, 2.64371) (0.62926, 2.64365) (0.63126, 2.64359) (0.63327, 2.64352) (0.63527, 2.64346) (0.63727, 2.64339) (0.63928, 2.64332) (0.64128, 2.64326) (0.64329, 2.64319) (0.64529, 2.64312) (0.64729, 2.64305) (0.64930, 2.64297) (0.65130, 2.64290) (0.65331, 2.64283) (0.65531, 2.64275) (0.65731, 2.64268) (0.65932, 2.64260) (0.66132, 2.64253) (0.66333, 2.64245) (0.66533, 2.64237) (0.66733, 2.64229) (0.66934, 2.64221) (0.67134, 2.64213) (0.67335, 2.64205) (0.67535, 2.64197) (0.67735, 2.64189) (0.67936, 2.64180) (0.68136, 2.64172) (0.68337, 2.64163) (0.68537, 2.64155) (0.68737, 2.64146) (0.68938, 2.64137) (0.69138, 2.64128) (0.69339, 2.64119) (0.69539, 2.64110) (0.69739, 2.64101) (0.69940, 2.64092) (0.70140, 2.64082) (0.70341, 2.64073) (0.70541, 2.64064) (0.70741, 2.64054) (0.70942, 2.64044) (0.71142, 2.64035) (0.71343, 2.64025) (0.71543, 2.64015) (0.71743, 2.64005) (0.71944, 2.63995) (0.72144, 2.63985) (0.72345, 2.63975) (0.72545, 2.63965) (0.72745, 2.63954) (0.72946, 2.63944) (0.73146, 2.63933) (0.73347, 2.63923) (0.73547, 2.63912) (0.73747, 2.63901) (0.73948, 2.63890) (0.74148, 2.63879) (0.74349, 2.63868) (0.74549, 2.63857) (0.74749, 2.63846) (0.74950, 2.63835) (0.75150, 2.63824) (0.75351, 2.63812) (0.75551, 2.63801) (0.75752, 2.63789) (0.75952, 2.63778) (0.76152, 2.63766) (0.76353, 2.63754) (0.76553, 2.63742) (0.76754, 2.63730) (0.76954, 2.63718) (0.77154, 2.63706) (0.77355, 2.63694) (0.77555, 2.63681) (0.77756, 2.63669) (0.77956, 2.63657) (0.78156, 2.63644) (0.78357, 2.63631) (0.78557, 2.63619) (0.78758, 2.63606) (0.78958, 2.63593) (0.79158, 2.63580) (0.79359, 2.63567) (0.79559, 2.63554) (0.79760, 2.63541) (0.79960, 2.63528) (0.80160, 2.63514) (0.80361, 2.63501) (0.80561, 2.63487) (0.80762, 2.63474) (0.80962, 2.63460) (0.81162, 2.63446) (0.81363, 2.63433) (0.81563, 2.63419) (0.81764, 2.63405) (0.81964, 2.63391) (0.82164, 2.63377) (0.82365, 2.63362) (0.82565, 2.63348) (0.82766, 2.63334) (0.82966, 2.63319) (0.83166, 2.63305) (0.83367, 2.63290) (0.83567, 2.63275) (0.83768, 2.63261) (0.83968, 2.63246) (0.84168, 2.63231) (0.84369, 2.63216) (0.84569, 2.63201) (0.84770, 2.63186) (0.84970, 2.63170) (0.85170, 2.63155) (0.85371, 2.63140) (0.85571, 2.63124) (0.85772, 2.63109) (0.85972, 2.63093) (0.86172, 2.63077) (0.86373, 2.63062) (0.86573, 2.63046) (0.86774, 2.63030) (0.86974, 2.63014) (0.87174, 2.62998) (0.87375, 2.62981) (0.87575, 2.62965) (0.87776, 2.62949) (0.87976, 2.62932) (0.88176, 2.62916) (0.88377, 2.62899) (0.88577, 2.62883) (0.88778, 2.62866) (0.88978, 2.62849) (0.89178, 2.62832) (0.89379, 2.62815) (0.89579, 2.62798) (0.89780, 2.62781) (0.89980, 2.62764) (0.90180, 2.62746) (0.90381, 2.62729) (0.90581, 2.62712) (0.90782, 2.62694) (0.90982, 2.62677) (0.91182, 2.62659) (0.91383, 2.62641) (0.91583, 2.62623) (0.91784, 2.62605) (0.91984, 2.62587) (0.92184, 2.62569) (0.92385, 2.62551) (0.92585, 2.62533) (0.92786, 2.62515) (0.92986, 2.62496) (0.93186, 2.62478) (0.93387, 2.62459) (0.93587, 2.62441) (0.93788, 2.62422) (0.93988, 2.62403) (0.94188, 2.62384) (0.94389, 2.62365) (0.94589, 2.62346) (0.94790, 2.62327) (0.94990, 2.62308) (0.95190, 2.62289) (0.95391, 2.62269) (0.95591, 2.62250) (0.95792, 2.62231) (0.95992, 2.62211) (0.96192, 2.62191) (0.96393, 2.62172) (0.96593, 2.62152) (0.96794, 2.62132) (0.96994, 2.62112) (0.97194, 2.62092) (0.97395, 2.62072) (0.97595, 2.62052) (0.97796, 2.62031) (0.97996, 2.62011) (0.98196, 2.61991) (0.98397, 2.61970) (0.98597, 2.61950) (0.98798, 2.61929) (0.98998, 2.61908) (0.99198, 2.61887) (0.99399, 2.61867) (0.99599, 2.61846) (0.99800, 2.61825) (1.00000, 2.61803) 
        };
    \end{axis}
\end{tikzpicture}
\end{example}

Despite of the simplicity of the statement, our proof of Theorem \ref{T:1} may seem rather indirect. In short, we will express the operator $\boldsymbol{m}$ in \eqref{E:T1} 
via a many-to-one formula. The latter involves  a multiplicative functional of a Markov chain which bears a simple connexion to the Perron--Frobenius eigen-elements of the mean reproduction matrix $m$. A key step is to establish unique ergodicity of this Markov chain, as then bounds for multiplicative functional can be deduced from Birkhoff's ergodic theorem. We stress that, since individuals in the population model keep full memory of the types of their ancestors, many  classical tools including irreducibility, strong Feller property, stability and contraction in total variation, compactness of operators, etc. (see e.g. the recent survey \cite{DMHJ} and also \cite{DMM}, and references therein) are not available. Nonetheless, we will see that this difficulty can be circumvented with 
the idea of asymptotic strong Feller property, which was developed notably in \cite{HM, HMS}.

The rest of this work is organized as follows. Section 2 essentially consists of notational preliminaries and a few simple observations about the Feller property and the spectral radius. Theorem \ref{T:1} is proved in Section 3, 
taking for granted  the unique ergodicity of the size-biased chain, which in turn is established in Section 4 by coupling. 

\section{Some notation and first observations}
Recall that  $S^\N$ denotes the space   of infinite sequences of types, and 
for every type $t\in S$ and memory $\s=(s_0, s_1, \ldots)$ in $S^\N$, we write $t\s=(t,s_0, s_1, \ldots)$ for the memory of a child with type $t$
of a parent with memory $\s$. Note that, in comparison with the Ulam-Harris-Neveu labeling for genealogies, concatenation is performed is the opposite order.
For every $\s, \s'\in S^\N$, we denote by $\s\wedge \s'$ the longest common prefix to $\s$  and $\s'$ (if any) and  $|\s\wedge \s'|$ for its length,
with the agreement that $|\s\wedge \s'|=0$ when $s_0\neq s'_0$ and that 
 $|\s\wedge \s'|=\infty$ if and only if $\s=\s'$. 
 
Recall that $T$ is the integrable random variable with law $\tau$ that selects the generation of the forebear whose type is activated for procreation. We write
 $$a_k=\P(T\geq k)=\sum_{j\geq k} \tau(j), \qquad k\in \N,$$
 and  also $a_{\infty}=0$.
So $(a_k: k\geq 0)$ is a non-increasing sequence in $[0,1]$ with $a_0=1< a_1$ and $\sum_{k=0}^{\infty} a_k<\infty$.
 We  equip the memory space $S^\N$ with the pseudo-distance
$$\dd(\s,\s')=a_{|\s\wedge \s'|},  \qquad \s, \s'\in S^\N.$$
More precisely, $\dd$ is a distance if and only if $a_k>0$ for all $k\geq 1$, that is, if and only if  the support of $\tau$ is unbounded,
and then we write $\Ss=S^\N$.
In the opposite case when $\max \mathrm{Supp}(\tau)=\ell<\infty$, we rather define $\Ss$ as the set of equivalence classes of $S^\N$ for $\dd$.
Then $\Ss$ can be identified as $S^\ell$; note also that the concatenation operation reads
$t(s_0, \ldots, s_{\ell-1})=(t, s_0, \ldots, s_{\ell-2})$.
In both cases, 
convergence of a sequence in $(\Ss,\dd)$ is equivalent to pointwise convergence, and 
 $(\Ss,\dd)$ is  clearly a compact metric space.

 We write $\mathcal B(\Ss)$ for the Banach algebra  of bounded measurable functions $f: \Ss\to \R$ equipped with the supremum norm, $\mathcal C(\Ss)$ for the sub-algebra of continuous functions, and  $ \Mm( \Ss)$ for the dual space of (finite) signed measures  on $\Ss$.  
 Recall that the operator $\boldsymbol{m}$ on $\mathcal B(\Ss)$ has been defined by \eqref{E:T1}; its adjoint $\boldsymbol{m}': \boldsymbol{\mu}\mapsto \boldsymbol{m}' \boldsymbol{\mu} = \boldsymbol{\mu} \boldsymbol{m}$  
is   a positive linear operator   on $ \Mm( \Ss)$ such that
$$\boldsymbol{\mu} \boldsymbol{m}(f)=\boldsymbol{\mu}( \boldsymbol{m} f )=\int_{\Ss} \boldsymbol{m}f(\s) \boldsymbol{\mu}(\dd \s),\qquad f\in \mathcal B(\Ss),$$
where positive means that $  \boldsymbol{\mu} \boldsymbol{m}$ is a positive measure  whenever $\boldsymbol{\mu}$ is a positive measure.
It may be interesting to make the following observation regarding the Feller property.

  \begin{proposition} \label{P:Feller}
The Feller property always holds, 
$$\boldsymbol{m}(\mathcal C(\Ss)) \subset \mathcal C( \Ss).$$
However,  the strong Feller property fails whenever  $\mathrm{Supp}(\tau)$ is unbounded, 
$$\boldsymbol{m}(\mathcal B( \Ss)) \not\subset \mathcal C( \Ss).$$ 
\end{proposition} 
\begin{proof}
The Feller property follows readily from \eqref{E:T1} and the  elementary inequality
$$\sum_{j=0}^{\infty} \tau(j) \indset{s_j\neq s'_j} \leq \dd(\s,\s').$$

The failure of
the strong Feller property when $\mathrm{Supp}(\tau)$ is unbounded stems  from the fact that the memory of parents can  be fully recovered from that of their children. 
More precisely, take any $\s\in \Ss$. Since $a_k>0$ for all $k\geq 1$, we may form a sequence $(\s(k))_{k\geq 1}$ in $\Ss$ such that $|\s\wedge \s(k)|=k$. So $\s(k)\neq \s$ and 
$\s(k)$ converges to $\s$ as $k\to \infty$.
 Consider  the subspace $\Ss'$ of memory sequences of the form $\s'=t\s$ for some $t\in S$ and take $f=\indset{\Ss'}$.
We have plainly $\boldsymbol{m}f(\s)=1$ whereas $\boldsymbol{m}f(\s(k))=0$ for each $k\geq 1$. 
\end{proof}

\begin{remark}
The argument justifying why strong Feller property fails also shows that $\boldsymbol{m}$ is not  a compact operator on $\mathcal C(\Ss)$
when $\mathrm{Supp}(\tau)$ is unbounded.
\end{remark}

In our population model with memory, the probability that an individual with memory $\s$ has $n_t$ children with type $t$ for every $t\in S$ (and no other children)
 equals 
 \begin{equation} \label{E:reprodlaw}
 \boldsymbol{\pi}(\s,n)=
 \sum_{j=0}^{\infty} \tau(j) \pi(s_j,n), \qquad n=(n_t)_{t\in S}\in \N^S.
 \end{equation}
We represent the population at generation $k$ as a point process $\Zz_k$ in $\Ss$, and  $\Zz=(\Zz_k: k\geq 0)$ is a branching process on $\Ss$ with reproduction law induced by $\boldsymbol{\pi}$. We write $\Pp(\Ss)$ for the space of probability measures on $\Ss$, and for every $\boldsymbol{\mu} \in \Pp(\Ss)$,  $\P_{\boldsymbol{\mu}}$ for the law of  $\Zz$ started from a single individual with random memory distributed according to $\boldsymbol{\mu}$.
 From \eqref{E:meanm}, \eqref{E:T1} and \eqref{E:reprodlaw}, the intensity measure of $\Zz_1$ under $\P_{\boldsymbol{\mu}}$ is $\boldsymbol{\mu} \boldsymbol{m}$.
 By iteration, we have for every $k\geq 0$ 
  $$\E_{\boldsymbol{\mu}}(\Zz_k(f))= \boldsymbol{\mu} \boldsymbol{m}^kf,$$
where we use the notation $\Zz_k(f)$ for the integral of $f\in \mathcal B( \Ss)$ with respect to the point process $\Zz_k$.

We write $\|\cdot\|$ for the operator norm on $\mathcal B( \Ss)$, so that 
$$\| \boldsymbol{m}^k\| = \sup_{\boldsymbol{\mu}\in \Pp(\Ss)} \ \boldsymbol{\mu} \boldsymbol{m}^k\mathbf{1}.$$
By Gelfand's formula, the spectral radius $\boldsymbol{r}$ of $\boldsymbol{m}$ can be expressed as
$$\boldsymbol{r}= \lim_{k\to \infty} \| \boldsymbol{m}^k\| ^{1/k},$$
and in particular, for any $c>\boldsymbol{r}$, we have
  $$\lim_{k\to \infty} c^{-k}\E_{\boldsymbol{\mu}}(\Zz_k(f))=0$$
  for any $\boldsymbol{\mu}\in \Pp(\Ss)$ and $f\in \mathcal{B}(\Ss)$. 
In the converse direction, we 
 point at the following weaker version of the Perron--Frobenius theorem.  
 
 \begin{proposition}\label{P:KRdual}  The spectral radius $\boldsymbol{r}$ is a positive eigenvalue of the adjoint  operator $\boldsymbol{m}'$.
 Namely $\boldsymbol{r}>0$ and there exists a probability measure $ \boldsymbol{\rho}$ on $ \Ss$  such that $ \boldsymbol{\rho} \boldsymbol{m} = \boldsymbol{r}  \boldsymbol{\rho}$. 
 \end{proposition}
  As a consequence, we have  that for any $k\geq 0$,
 $$\E_{\boldsymbol{\rho}}(\Zz_k(\mathbf{1}))= \boldsymbol{r}^k,$$
 which confirms that the spectral radius $\boldsymbol{r}$ indeed determines the growth rate on average of the population. 

 \begin{proof} 
 The case when $\mathrm{Supp}(\tau)$ is bounded is plain from the Perron--Frobenius theorem, since then $\Ss$ a finite set. Hence we focus on the case when $\mathrm{Supp}(\tau)$ is unbounded; the idea of the proof is to establish first the claim for a truncated version of $\boldsymbol{m}$ and then conclude by taking a limit. 

So fix any $\ell \geq 1$ and introduce the operator $\boldsymbol{m}_{\ell}$ on $\mathcal{B}(\Ss)$ given by 
$$
   \boldsymbol{m}_{\ell}f(\s)= \sum_{j=0}^{\ell-1}\tau(j) \sum_{t\in S} m(s_j,t) f(t\s), \qquad \s \in \Ss.
$$
Note that for any $k\geq \ell$, the subspace $\R^{S^k}\subset \mathcal{B}(\Ss)$ of functions $f(\s)=f(s_0, \ldots, s_{k-1})$ that only depend on the first $k$ elements of the memory $\s$
is invariant for $ \boldsymbol{m}_{\ell}$. By restriction, we can also view $\boldsymbol{m}_{\ell}$ as a positive linear operator on the finite dimensional space $\R^{S^k}$.
We stress that by Gelfand's formula, this does not affect 
 the spectral radius $\boldsymbol{r}_{\ell}$ of $\boldsymbol{m}_{\ell}$, that is the spectral radius  remains the same, no matter whether $\boldsymbol{m}_{\ell}$ is viewed as an operator on $\mathcal{B}(\Ss)$ or on $\R^{S^k}$, 
 since  $\boldsymbol{m}_{\ell}\mathbf{1}\in \R^{S^k}$ for every $k\geq \ell$.
Plainly, there are the bounds
$$0<  \min_{s\in S} \sum_{t\in S} m(s,t) \leq  \boldsymbol{r}_{\ell} \leq \  \max_{s\in S} \sum_{t\in S} m(s,t) <\infty,$$ and by the version of the Perron--Frobenius theorem for (possibly reducible) matrices with nonnegative entries, for any $k\geq \ell$,
there exists ${\boldsymbol{\rho}}_{k}\in \Pp(S^{k})$ such that 
\begin{equation} \label{E:eigenk}{\boldsymbol{\rho}}_{k} \boldsymbol{m}_{\ell}= \boldsymbol{r}_{\ell}{\boldsymbol{\rho}}_{k}.
\end{equation}

Our next goal is to extend  \eqref{E:eigenk} to a solution of the eigen-problem when $\boldsymbol{m}_{\ell}$ is viewed as an operator on $\mathcal{B}(\Ss)$ rather than on $\R^{S^k}$.
In this direction, observe that if ${\boldsymbol{\rho}}'_{k}$ is the image of ${\boldsymbol{\rho}}_{k+1}$ by the projection $(s_0, \ldots, s_k) \mapsto (s_0, \ldots, s_{k-1})$,
then ${\boldsymbol{\rho}}'_k$ also solves \eqref{E:eigenk}. By compactness of $\Pp(S^{k})$ and an argument of diagonal extraction, 
we may now further assume that the eigen-laws ${\boldsymbol{\rho}}_{k}$ have been constructed such that ${\boldsymbol{\rho}}_{k}$ is the image of ${\boldsymbol{\rho}}_{k+1}$ by the projection $(s_0, \ldots, s_k) \mapsto (s_0, \ldots, s_{k-1})$ for any $k\geq \ell$.
By Kolmogorov's extension theorem, there is a probability measure $ \boldsymbol{\rho}_{\ell}$ on $S^\N=\Ss$  whose image  
by the projection $\s\mapsto (s_0, \ldots, s_{k-1})$ is ${\boldsymbol{\rho}}_{k}$ for any $k\geq \ell$. 
 Plainly, now seeing  $\boldsymbol{m}_{\ell}$ as an operator 
on $\mathcal{B}(\Ss)$, we have $ \boldsymbol{\rho}_{\ell} \boldsymbol{m}_{\ell}= \boldsymbol{r}_{\ell} \boldsymbol{\rho}_{\ell}$. In other words, $ \boldsymbol{\rho}_{\ell} $ is an eigen-law for the eigenvalue $\boldsymbol{r}_{\ell}$ of  the adjoint of the operator $\boldsymbol{m}_{\ell}$ on $\mathcal{B}(\Ss)$. 

To complete the proof, we let $\ell\to \infty$.  Observe that
$$ \|  \boldsymbol{m}_{\ell}-  \boldsymbol{m} \| \leq a_\ell \times \max_{s\in S} \sum_{t\in S} m(s,t);$$
 it follows from Gelfand's formula  that  $\boldsymbol{r}_{\ell}$ converges to the spectral radius $ \boldsymbol{r}$ of $ \boldsymbol{m}$. 
On the other hand, since $\Ss$ is compact, we can find some  $ \boldsymbol{\rho} \in \Pp(\Ss)$ and extract a subsequence such that $ \boldsymbol{\rho}_{\ell}\to  \boldsymbol{\rho}$ weakly
and then $ \boldsymbol{\rho} \boldsymbol{m} = \boldsymbol{r}  \boldsymbol{\rho}$.
 \end{proof}

Proposition \ref{P:KRdual} is only a first step towards what one may expect, namely the full extension of the Perron--Frobenius theorem which would include the existence of a strictly positive  eigen-function $\boldsymbol{h}\in \mathcal{C}(\Ss)$ with $ \boldsymbol{\m} \boldsymbol{h} = \boldsymbol{r}  \boldsymbol{h}$, and stability and contraction theorems; see, e.g. \cite[Theorem 21.1]{BW} or \cite[Section 4]{DMHJ}. Of course, this plainly holds when $\mathrm{Supp}(\tau)$ is bounded, but remained in general elusive for the author when $\mathrm{Supp}(\tau)$ is unbounded. See nonetheless Remark \ref{R:CSV} in the next section.

\section{A biased Markov chain}
This section aims to establish Theorem \ref{T:1} via an identity known as a many-to-one formula,
which relates  powers of the  operator $\boldsymbol{m}$ in \eqref{E:T1} to a multiplicative functional 
 a so-called biased Markov chain. There is a great variety of such formulas, and our choice is guided by the spectral properties of the model without memory. 
 The existence of a unique invariant law for the biased chain lies at the heart of the proof; this will be established in the next section.

It is convenient to record for each individual, not only its memory $\s\in \Ss$, but also the type $s=s_T$ of its forebear that has been activated for its procreation. We refer to the pair $(s,\s)$ as the symbol of an individual, and write $\overline \Ss=S\times \Ss$ for the space of symbols. 
Note  that the first element of the symbol
has no impact on the reproduction of this individual.

Recall that $r>0$ denotes the spectral radius of the mean reproduction matrix $m$, which is irreducible and aperiodic by assumption, and that the Perron--Frobenius theorem
ensures the existence of a unique solution $(\rho, h)$ of the eigen-problem
$$ \rho m= r \rho\quad \text{and} \quad mh = rh,$$
where $\rho$ is a  probability measure on $S$ and  $h: S\to (0,\infty)$ a function 
such that
$$\rho(h)=\sum_{s\in S} \rho(s) h(s)=1.$$
This
incites us to
 introduce
  the normalized  transition matrix $\overline m$ in $\R_+^{S\times S} $ with entries
\begin{equation} \label{E:mtilde}
\overline m(s,t)=\frac{m(s,t)h(t)}{rh(s)}, \qquad s,t\in S,
\end{equation}
and then, mimicking \eqref{E:T1},   the  Markovian transition kernel 
\begin{equation}\label{E:barqbis}
\overline Qf(s,\s) =  \sum_{j=0}^{\infty} \tau(j)  \sum_{t\in S}\overline m(s_j,t)f(s_j,t\s),\qquad f\in \mathcal{B}(\overline \Ss),
\end{equation}
where $\mathcal{B}(\overline \Ss)$ stands for the space of bounded measurable functions on $\overline \Ss= S\times \Ss$.
In words, the evolution consists in  first activating a type $s_j$ at random from the memory using the law $\tau$, and then following the normalized transition 
$\overline{m}(s_j, \cdot)$. 
We stress that the right-hand side in \eqref{E:barqbis} purposely does not depend on the initial type $s$.

We write $\overline \X=((Y_k,\X_k): k\geq 0)$ for the Markov chain in $\overline \Ss$ with transition kernel $\overline Q$ and call it the biased chain.
Recall that the initial value $Y_0=s$ plays no role for the rest of the evolution, and for simplicity, we may sometimes omit it from notation when it is not needed.
For instance, we  interchangeably write $\overline \P_{\s}$ or $\overline \P_{\overline{\s}}$  for the law of the chain started  from $\X_0=\s$ 
or from $(Y_0, \X_0)=\overline \s=(s,\s)$, keeping in mind that the first component  $s$
of the symbol $\overline \s$ is used for $Y_0$ only, and is irrelevant for the other variables appearing in the process. 
We also write $X_k$ for  the first element of the sequence $\X_k$, so that  $(X_k, X_{k-1}, \ldots, X_0)$ is the prefix with length $k+1$ of $\X_k$.
We may think of $X_k$ as the type of a distinguished individual at time $k$; beware also that $((Y_k,X_k): k\geq 0)$ is not Markovian. 
In this framework, the many-to-one formula 
\begin{equation}\label{E:many21}
\boldsymbol{m}^kf(\s)= r^k\overline \E_{\s}\left(f(\X_k) \times  \prod_{j=1}^{k} \frac{h(Y_j)}{h(X_j)}\right),
\qquad f\in \mathcal{B}(\Ss),
\end{equation}  follows easily from \eqref{E:T1},  \eqref{E:mtilde} and \eqref{E:barqbis} by induction on $k\geq 0$.
We first point at  elementary bounds for the spectral radii
\begin{equation} \label{E:Harn}
 \frac{\min h}{\max h}  \leq \frac{\boldsymbol{r}}{r} \leq \frac{\max h}{\min h}.
 \end{equation} 
\begin{example} \label{E:bal} If the mean reproduction matrix $m$ is balanced, in the sense that the mean offspring $\sum_{t\in S} m(s,t)$ when type $s$ is activated does not depend on $s$,
then the latter quantity coincides with  the spectral radius $r$ of $m$ and $h\equiv 1$. We deduce that the spectral radius of $\boldsymbol{m}$ is $\boldsymbol{r}=r$. 
\end{example}

The many-to-one formula suggests that one should gain insight on $\boldsymbol{m}^k$ when $k$ is large by analyzing the asymptotic behavior of the Markov chain $\overline \X$. 
  \begin{proposition} \label{P:uniferg}
 The  kernel $\overline Q$ is uniquely ergodic, i.e.  there exists a unique law $\overline {\boldsymbol{\sigma}}\in \Pp(\overline \Ss)$ that solves $\overline {\boldsymbol{\sigma}}= \overline {\boldsymbol{\sigma}} \overline Q$.  
  \end{proposition}

Proposition \ref{P:uniferg} is intuitive when the support of $\tau$ is bounded, as the state space $\overline{\Ss}$ becomes finite (although irreducibility may not always hold). However, the scenario becomes less clear when the support is unbounded; the proof of this case is deferred to the next section.
 
 The next ingredient needed for the proof of Theorem \ref{T:1}
 is the observation that the one-dimensional marginals of the invariant law $\overline {\boldsymbol{\sigma}}$ are all the same.
 
 \begin{lemma} \label{L:marg} For every $t\in S$ and $j\geq 0$, there are the identities
 $$\overline {\boldsymbol{\sigma}}\left( \left \{(s,\s)\in \overline \Ss: s=t\right \}\right) = 
 \overline {\boldsymbol{\sigma}}\left( \left \{(s,\s)\in \overline \Ss: s_j=t \right \}\right) = \rho(t) h(t). $$
 More generally, we have for the two-dimensional marginal that
 $$\overline {\boldsymbol{\sigma}}\left( \left \{(s,\s)\in \overline \Ss: s=t \text{ and } s_0=u\right \}\right) = \rho(t)m(t,u) h(u)/r, \qquad t,u\in S.$$
 \end{lemma}
 \begin{proof} Since  $\overline {\boldsymbol{\sigma}}= \overline {\boldsymbol{\sigma}}\overline Q^j$  for any  $j\geq 0$,
 the marginal law $\sigma_j$ of $s_j$ under $ \overline {\boldsymbol{\sigma}}$ does not depend on $j$, and we simply denote it by $\sigma$.
 It then follows from \eqref{E:barqbis} that for any $t\in S$,
 we can express $\sigma_1(t)=\sigma(t)$ in the form
 \begin{align*}\int_{\overline \Ss}  \overline {\boldsymbol{\sigma}} (\dd s, \dd \s) \left(   \sum_{j=0}^{\infty} \tau(j)  \overline m(s_j,t) \right)
 &=   \sum_{j=0}^{\infty} \tau(j) \sum_{s'\in S} \sigma_j(s')  \overline m(s',t)\\
 & = \sum_{s'\in S} \sigma(s') \overline m(s',t),
 \end{align*}
 where for the last equality, we used that $\sigma_j=\sigma$ and that $\tau$ is a probability measure.
We get from \eqref{E:mtilde}
$$r\frac{\sigma(t)}{h(t)}= \sum_{s'\in S}  \frac{\sigma(s')}{h(s')}m(s',t).$$
Since $\rho$ is the unique probability measure on $S$ with $\rho m = r \rho$, 
we conclude that $\sigma(t)=h(t)\rho(t)$.
 
Similarly, we have from \eqref{E:barqbis} that for any $s'\in S$, there are the identities
$$
\overline {\boldsymbol{\sigma}}\left( \left \{(s,\s)\in \overline \Ss: s=t\right \}\right) =  \sum_{j=0}^{\infty} \tau(j)  \sum_{t\in S}  \sigma_j(t) 
= \sigma(t). 
$$
The calculation for the two-dimensional marginal is also the same.
 \end{proof}
 
 Let us now establish Theorem \ref{T:1}, taking  Proposition \ref{P:uniferg} for granted.
 
 \begin{proof}[Proof of Theorem \ref{T:1}]
 It suffices to check that 
 \begin{equation} \label{E:proofT1}
\liminf_{k\to \infty}\frac{1}{k}  \sum_{j=1}^{k} (\log h(Y_j)- \log h(X_j))\geq 0, \qquad \overline{\P}_{\overline {\boldsymbol{\sigma}}}\text{-a.s.}
\end{equation}
Indeed, assuming \eqref{E:proofT1} for the time being, we rewrite it as
$$\liminf_{k\to \infty}\left(  \prod_{j=1}^{k} \frac{h(Y_j)}{h(X_j)}\right)^{1/k}\geq 1, \qquad \overline{\P}_{\overline {\boldsymbol{\sigma}}}\text{-a.s.}
$$
Combining Jensen's inequality and Fatou's lemma, we arrive at
$$\liminf_{k\to \infty} \left(\overline \E_{\overline {\boldsymbol{\sigma}}}\left(  \prod_{j=1}^{k} \frac{h(Y_j)}{h(X_j)}\right)\right)^{1/k}\geq 1.$$
Next, writing ${\boldsymbol{\sigma}}$ for image of 
$\overline {\boldsymbol{\sigma}}$ by the second projection, $(s,\s)\mapsto \s$,  from $\overline \Ss$ to $\Ss$, 
 the many-to-one formula \eqref{E:many21}  
 yields 
$$ \boldsymbol{\sigma} \boldsymbol{m}^k \mathbf{1}= r^k\overline \E_{\overline {\boldsymbol{\sigma}}}\left(  \prod_{j=1}^{k} \frac{h(Y_j)}{h(X_j)}\right),$$
and therefore
$$\boldsymbol{r}= \lim_{k\to \infty}  \sup_{\boldsymbol{\mu}\in \Pp(\Ss)} \left(\boldsymbol{\mu} \boldsymbol{m}^k\mathbf{1}\right)^{1/k}\geq r.$$

It remains to prove  \eqref{E:proofT1}. We know from Proposition \ref{P:uniferg} that the law $ \overline {\boldsymbol{\sigma}}$
is $\overline Q$-invariant and ergodic; see, e.g. \cite[Proposition 4.30]{BH}. The law $\overline{\P}_{\overline {\boldsymbol{\sigma}}}$ of the
Markov chain $\overline \X$ is shift-invariant and ergodic; see, e.g. \cite[Proposition 4.49]{BH}. 
By Birkhoff ergodic theorem, see, e.g. \cite[Theorem 4.42]{BH}), 
  $\overline{\P}_{\overline {\boldsymbol{\sigma}}}$-a.s. there is the convergence
$$\lim_{k\to \infty}\frac{1}{k}  \sum_{j=1}^{k} (\log h(Y_j)- \log h(X_j))= \int_{\overline \Ss} (\log h(s)-\log h(s_0)) \overline {\boldsymbol{\sigma}}(\dd s , \dd \s).$$
We now  see from Lemma \ref{L:marg} that the right-hand side equals $0$, which completes the proof. 
 \end{proof}

 We now conclude this section with some further comments about the eigen-problem for $\boldsymbol{m}$ at the light of the many-to-one formula and of  the recent work by Champagnat, Strickler and Villemonais \cite{CSV}. 

\begin{remark}\label{R:CSV} 
Consider $c>0$ sufficiently small such that $h(s) \geq c h(t)$ for all $s,t\in S$; we can re-write the many-to-one formula \eqref{E:many21} in the form
$$
\overline \E_{\s}\left(f(\X_k) \times  \prod_{j=1}^{k} p(\overline \X_k)\right)= c^{k} r^{-k} \boldsymbol{m}^kf(\s),
$$  
where $p(s,\s)= c h(s)/h(s_0)\in (0,1)$.
The left-hand side is interpretated as a penalization of the Markov chain $\overline \X$.
The main purpose of the authors in \cite{CSV} is to 
study for fairly general penalized Markov processes the convergence to a quasi-stationary distribution.
They have obtained a simple criterion for the uniform convergence  in Wasserstein distance
 of conditional distributions given survival; see \cite[Theorem 2.1]{CSV}, and also shown that under the same requirement, 
 the principal eigen-problem has a solution; see \cite[Theorem 2.3]{CSV}. 
 Unfortunately, their key assumption, namely Assumption (A), seems to be quite challenging to verify in our setting;
 the same holds for the weaker Harnack type inequality (H) there. 
\end{remark}
 
 \section{An asymptotic coupling}
 The purpose of this section is to establish
Proposition \ref{P:uniferg}. For this, we will rely on
coupling,  a well-know technique introduced by Doeblin to show convergence to equilibrium for positive recurrent Markov chains on countable state spaces, and then extended by Harris to Markov chains on general state spaces; see e.g.  \cite{BH, BW0, MeynT}. The original idea is to construct simultaneously  two trajectories of the Markov chain with arbitrary starting states that 
 eventually meet at the same point and merge thereafter.
Such exact couplings are impossible in our setting, since the chain keeps a full memory of its past,
in the sense that for any $0\leq k \leq k'$,  $\X_k$ is recovered as the $(k'-k)$-th suffix of $\X_{k '}$ (in particular $\overline \X$ is not  a Harris chain).
Therefore  we shall rather couple the trajectories such that after some random finite time,
 they keep getting closer and closer. Even though the two trajectories never meet, it is sufficient for our purposes. This idea of
 asymptotic coupling has been developed notably by Hairer-Mattingly \& Scheutzow \cite{HMS}; see also the references therein for earlier works.

   Technically  in our setting, a coupling consists of a Markovian transition kernel $ C$ on $\overline \Ss \times \overline \Ss$ which yields for any initial condition $(\overline \s, \overline \s')\in  \overline \Ss \times  \overline \Ss$
 the law $ \C_{\overline \s, \overline \s'}$ of 
 a bivariate chain $(\overline \X,\overline \X')$  such that $\overline \X$ has the law $\overline \P_{\overline \s}$ and $\overline \X'$ the law $\overline \P_{\overline \s'}$. The efficiency of a coupling
is evaluated by comparing the pairs of types $(Y_k,X_k)$ and $(Y'_k,X'_k)$ (recall that $X_k$ denotes the first element of the sequence $\X_k$) for $k\geq 1$.
Given a realization of the coupling, we say that the $k$-th step is successful if $(Y_k,X_k)=(Y'_k,X'_k)$, and fails otherwise. 
We say that a coupling is ultimately successful if for any  pair $(\overline \s, \overline \s')$ of initial conditions, 
all but finitely many steps are successful, a.s.
The main purpose of this section is to check the following existence claim. 
 
 \begin{lemma} \label{L:upe} There exists an ultimately successful Markovian coupling. \end{lemma} 
 
 Let us briefly deduce Proposition \ref{P:uniferg}  from Lemma \ref{L:upe}.
 \begin{proof}[Proof of Proposition \ref{P:uniferg}]
 Just as in Proposition \ref{P:Feller}, we observe that $\overline Q$ is a Feller operator. Since $\overline \Ss$ is compact, the existence of an invariant law $\overline {\boldsymbol{\sigma}}$ follows classically from the Schauder fixed point theorem.

 To prove uniqueness, consider a second invariant law $\overline {\boldsymbol{\sigma}}'$.
  Given an ultimately successful coupling $C$,  introduce the mixture
 $$\C_{\overline {\boldsymbol{\sigma}}, \overline {\boldsymbol{\sigma}}'}= \int_{\overline \Ss\times \overline \Ss} \overline {\boldsymbol{\sigma}} (\dd \overline \s) \overline {\boldsymbol{\sigma}}'(\dd \overline \s')  \C_{\overline \s, \overline \s'}.$$
  So under $\C_{\overline {\boldsymbol{\sigma}}, \overline {\boldsymbol{\sigma}}'}$,  for every $k\geq 1$, the law of $\overline \X_k$  is $\overline {\boldsymbol{\sigma}}$ and  the law of $\overline \X'_k$  is $\overline {\boldsymbol{\sigma}}'$. 
Take any $f\in \mathcal{B}(\overline \Ss)$ such that  the function $f(s, \cdot): \Ss\to [0,1]$ is $1$-Lipschitz-continuous for every $s\in S$.
Then 
$$\overline {\boldsymbol{\sigma}}(f)  - \overline {\boldsymbol{\sigma}}' (f)  =  \E (f(\overline \X_k) -f(\overline \X'_k)) ,$$
where the expectation in the right-hand side is computed under
$\C_{\overline {\boldsymbol{\sigma}}, \overline {\boldsymbol{\sigma}}'}$.
Now it suffices to write
$$|f(\overline \X_k) -f(\overline \X'_k)| \leq \dd(\X_k, \X'_k) + \indset{Y_k\neq Y'_k},$$  
and deduce from Lemma \ref{L:upe} that as $k\to \infty$, the right-hand side converges to $0$ in probability.
We conclude by dominated convergence that $\overline {\boldsymbol{\sigma}}(f) = \overline {\boldsymbol{\sigma}}' (f)$
for all such functions $f$, and hence $\overline {\boldsymbol{\sigma}}= \overline {\boldsymbol{\sigma}}'$. 
  \end{proof}

 Our purpose now is to check Lemma \ref{L:upe} and we construct below  a bivariate Markovian transition kernel $C$ on $\overline \Ss\times \overline \Ss$.
   The idea is to use the common parts of the memories $\X_k$ and $\X'_k$ to activate the same type  and then let the two chains evolve in the same way for the next step. 
   We iterate as long as possible.
   
   \begin{proof}[Proof of  Lemma \ref{L:upe}]
 We only need to specify the coupling operator  $C$ on the product space $\mathcal B(\overline \Ss)\times \mathcal B(\overline \Ss)$,
 and our description focusses on functions $F: \overline \Ss^2\to \R$ of the type $F(\overline \s, \overline \s')=f(\overline \s)f'(\overline \s')$ for $f,f'\in \mathcal B(\overline \Ss)$
 and symbols given in the form $\overline\s=(s,\s)$ and $\overline \s'=(s',\s')$. We
 write $\ell=|\s\wedge \s'|$ for the length of the longest common prefix of the memories  
  and define
\begin{align*}
C F(\overline \s, \overline \s') =& \sum_{j=0}^{\ell-1} \tau(j) \sum_{t\in S} \overline m(s_j,t) f(s_j, t\s) f(s_j, t \s')\\
&+ \sum_{j=\ell}^{\infty} \tau(j)  \left(
 \sum_{t\in S} \overline m(s_j,t) f(s_j, t\s)\right) \left(  \sum_{t\in S} \overline m(s'_j,t) f(s'_j, t\s')\right);
\end{align*}
we stress that the first sum in the right-hand side vanishes when $\ell=0$.

In words, we distinguish two phases of the coupling. The initialization phase is when $\ell=0$; then the two chains evolve independently one from the other.
The consolidation phase starts when $\ell\geq 1$; then we use the longest common prefix of the memory to activate the same type and further
pair their evolutions, but we also keep
the possibility of letting the two chains make different steps in order to get the right marginals in the end.
In short, we will argue that the coupling always exit the initialization phase after some random number of steps with bounded expectation,
and that once the coupling entered  the consolidation phase, the probability that it remains there forever is bounded away from $0$.
Even though the coupling may leave the consolidation phase and enter again the initialization phase, the number of such transitions is stochastically bounded 
by a geometric variable; this enables us to conclude that $C$ is indeed ultimately successful.

To start with, we immediately check from \eqref{E:barqbis} that under $\C_{\overline \s, \overline \s'}$, $\overline \X$ has indeed the law $\overline \P_{\overline \s}$ and $\overline \X'$ 
the law $\overline \P_{\overline \s'}$. 
We next consider the first entrance time into the consolidation phase,
$$\gamma=\inf\{k\geq 0: X_k=X'_k\}.$$
The transition matrix $\overline m$ in \eqref{E:mtilde} is irreducible and aperiodic, since this is the case for $m$ and $h>0$; there are  $\epsilon >0$ and $k_0\geq 1$ such that
$$\inf_{s,t\in S}\overline m^{k_0}(s,t) > \epsilon.$$
Recalling that $\tau(0)>0$, we now see by focussing on the event when the $k_0$ first steps of $\overline \X$ and $\overline \X'$ do not use their memories (in the sense that the activated type for the transition is simply the type of the parent)  that 
$$\C_{\overline \s, \overline \s'}(\gamma \leq  k_0) \geq \tau(0)^{2k_0}\epsilon,\qquad \text{for any } \overline \s, \overline \s'\in \overline \Ss.$$
We deduce by iteration from the Markov property that there is the uniform bound
$$ \overline \E_{\overline \s, \overline \s'}(\gamma) \leq \frac{ k_0}{\tau(0)^{2k_0}\epsilon},$$
where the expectation above is taken under the coupling law $\C_{\overline \s, \overline \s'}$.

We turn our attention to the consolidation phase. 
Recall that 
$$a_k=\sum_{j=k}^{\infty} \tau(j) \quad \text{and} \quad \sum_{k\geq 1} a_k<\infty.$$
We deduce from the definition of the coupling and an iteration using the Markov property that for any $\overline \s, \overline \s'\in \overline \Ss$ with $|\s\wedge \s'|=\ell\geq 1$ and
 any $n\geq 1$, we have
$$\C_{\overline \s, \overline \s'}\left((Y_k, X_k)= (Y'_k, X'_k) \text{ for all } k=1, \ldots, n 
\right) \geq \prod_{k= \ell}^{\ell + n-1} (1- a_k).$$
We thus have the uniform lower bound in the consolidation phase
$$\C_{\overline \s, \overline \s'}\left((Y_k, X_k)= (Y'_k, X'_k) \text{ for all } k\geq 1
\right) \geq \prod_{k= 1}^{\infty} (1- a_k) >0.$$
Putting the pieces together, this confirms that $C$ is indeed ultimately successful.  
  \end{proof}
 
\bibliographystyle{plain}
\bibliography{biblio}

\end{document}